\documentclass[11pt]{amsart}

\usepackage[T1]{fontenc}
\usepackage[utf8]{inputenc}
\usepackage{amsmath,amssymb,amsthm}
\usepackage{enumitem}
\usepackage[a4paper,margin=3cm]{geometry}
\usepackage{placeins}
\usepackage{tikz}
\usepackage[hidelinks,pdfencoding=auto]{hyperref}

\newtheorem{theorem}{Theorem}[section]
\newtheorem{proposition}[theorem]{Proposition}
\newtheorem{corollary}[theorem]{Corollary}

\theoremstyle{definition}
\newtheorem{example}[theorem]{Example}

\title[A Note on ``Localization of Zeros of Polar Polynomials on the Unit Disk'']
{A Note on the Paper ``Localization of Zeros of Polar Polynomials on the Unit Disk''}

\author{R. \'Alvarez-Nodarse}
\address{Departamento de Análisis Matemático, Universidad de Sevilla}
\address{c/Tarfia s/n, 41012 Sevilla, Spain}
\email{ran@us.es}

\author{K. Castillo}
\address{CMUC, Department of Mathematics, University of Coimbra,
3000-143 Coimbra, Portugal}
\email{kenier@mat.uc.pt}

\date{29 July 2026}

\hypersetup{
  pdftitle={A Note on the Paper ``Localization of Zeros of Polar
  Polynomials on the Unit Disk''},
  pdfauthor={R. Alvarez-Nodarse and K. Castillo},
  pdfsubject={A mathematical correction and a sharp Grace--Szego
  localisation theorem},
  pdfkeywords={zeros of polynomials, Grace--Szego convolution,
  orthogonal polynomials on the unit circle}
}

\begin{document}

\begin{abstract}
We show that none of the principal claims of Costas-Santos and Rehouma,
\emph{Rocky Mountain J. Math.} \textbf{54} (2024), 995--1004, survives as
a new and valid result.  The main localisation theorem is false as stated;
under the missing hypothesis of orthogonality on the unit circle, it is a
classical exercise recorded by Marden and by Borwein and Erdélyi, while
the published proof omits its decisive finite-degree estimate.  The
annular claim and Sendov application are unproved, and some reported
numerical values contradict the latter under its required hypothesis.
We establish instead a measure-independent Grace--Szeg\H{o} theorem with
sharp affine and optimal disc bounds, and characterise exactly when the
convolution method extends to a general factor.  This provides the
necessary correction to the published record.
\end{abstract}

\maketitle

\section{Introduction}

This note corrects the principal mathematical claims of \cite{CR24}.
Theorem~3.2 is false under its stated hypotheses.  When the intended
unit-circle orthogonality hypothesis is imposed, its conclusion reduces to
a classical textbook exercise, but the published proof still omits the
finite-degree estimate on which the stated radius depends.  The asserted
annular localisation is not established, and the Sendov discussion neither
verifies the hypothesis of the conjecture nor supports its conclusion.
Since these statements form part of the published literature, their
correction should likewise be placed on record.

Costas-Santos and Rehouma state the following result
\cite[Theorem~3.2]{CR24}.
\begin{quote}
\emph{\textbf{Theorem 3.2.}
Let $\mu$ be a finite measure defined on the Borelian $\sigma$-algebra
of $\mathbb C$ such that it contains an infinite number of points and let
$(L_n(z))$ be the system of monic orthogonal polynomials with respect to
$\mu$. Let $\xi$ be a fixed complex number and let $k$ be a positive
integer. All the zeros of $Q_{n;k}(z;\xi)$ are contained in the closed
disk
\[
 \overline{D\!\left(0,|\xi|+(k+1)(1+|\xi|)\right)}.
\]}
\end{quote}
Throughout,
\[
 D(c,r):=\{z\in\mathbb C:|z-c|<r\},\quad
 \overline D(c,r):=\{z\in\mathbb C:|z-c|\le r\},
\]
and \(Z(F)\) denotes the zero set of a polynomial \(F\), with
multiplicities retained whenever explicitly stated.
Using the rising factorial
\((a)_k:=a(a+1)\cdots(a+k-1)\), their Definition~1.3 introduces
$Q_{n;k}$ as the degree-$n$ polynomial solution of
\begin{equation}
 \frac{\mathrm d^k}{\mathrm dz^k}
 \bigl((z-\xi)^k Q_{n;k}(z;\xi)\bigr)
 =(n+1)_kL_n(z),
\end{equation}

There are two logically distinct defects in these hypotheses.

\begin{enumerate}[label=\textup{(\roman*)}]
\item Even if ``measure'' is read in its standard positive sense and
``contains infinitely many points'' is read as ``has infinite support'',
a finite Borel measure on $\mathbb C$ need not have finite moments of all
orders.  Thus finiteness of $\mu$ alone does not ensure that the inner
products required to construct $L_n$ are finite.  For a fixed degree
\(n\), a positive measure satisfying
\[
 \int_{\mathbb C}|z|^{2n}\,\mathrm d\mu(z)<\infty
\]
and having at least \(n+1\) support points is necessary and sufficient
for the usual Gram--Schmidt construction through degree \(n\).  For an
orthogonal sequence in every degree, the corresponding necessary and
sufficient conditions are that a positive Borel measure have infinite
support and satisfy
\[
 \int_{\mathbb C}|z|^{2m}\,\mathrm d\mu(z)<\infty,\quad m\ge0.
\]
Indeed, infinite support then ensures that
\[
 \int_{\mathbb C}|p(z)|^2\,\mathrm d\mu(z)>0,\quad
 p\in\mathbb C[z],\quad p\ne0,
\]
because a non-zero one-variable polynomial has only finitely many zeros.
Thus the monic orthogonal polynomials exist in every degree.  If the
finite measure is supported on \(\mathbb T\), all moments are finite
automatically; positivity and infinite support are then the essential
conditions for a full OPUC sequence.  If the words ``let \((L_n)\) be the
system'' are instead read as an independent existence assumption, the
missing moment and non-degeneracy requirements are implicit rather than
consequences of the stated assumptions on \(\mu\).

\item Even those corrected assumptions do not make $(L_n)$ a sequence of
orthogonal polynomials on the unit circle (OPUC).  For that, one must impose
\(\operatorname{supp}\mu\subset\mathbb T\), where
\(\mathbb T=\{z:|z|=1\}\).  The zero property used in the proof---that every
zero of $L_n$ lies in the open unit disc---is an OPUC theorem and therefore
cannot be invoked for an arbitrary measure on $\mathbb C$.
\end{enumerate}

The defect in \textup{(ii)} is not merely a gap in the published proof:
under the literal hypotheses the asserted disc is false.

\begin{example}[Counterexample to the theorem as stated]
\label{ex:translated-circle}
Fix \(A\in\mathbb C\) with \(|A|>1\), and let \(\mu_A\) be normalised arc
length on the translated circle \(A+\mathbb T\); equivalently, \(\mu_A\)
is the push-forward of normalised arc length on \(\mathbb T\) under
\(\zeta\mapsto A+\zeta\).  This is a finite positive compactly supported
Borel measure with infinite support and moments of every order.  Its monic
orthogonal polynomials are
\[
 L_n(z)=(z-A)^n,
\]
since, for \(m,n\ge0\),
\[
 \int_{A+\mathbb T}(z-A)^n
 \overline{(z-A)^m}\,\mathrm d\mu_A(z)
 =
 \int_{\mathbb T}\zeta^n\overline{\zeta^m}\,
 \frac{\lvert\mathrm d\zeta\rvert}{2\pi}
 =\delta_{mn}.
\]
Take \(n=1\), \(\xi=0\), and any \(k\ge1\).  The unique monic solution is
\[
 Q_{1;k}(z;0)=z-(k+1)A,
\]
because
\[
 \frac{\mathrm d^k}{\mathrm dz^k}
 \left(z^k\bigl(z-(k+1)A\bigr)\right)
 =(k+1)!(z-A)=(2)_kL_1(z).
\]
Its zero has modulus \((k+1)|A|>k+1\), whereas the disc asserted in
\cite[Theorem~3.2]{CR24}, with \(\xi=0\), has radius \(k+1\).
Thus the published theorem fails even for a positive, compactly supported
measure for which the complete orthogonal sequence exists.
\end{example}

The paper also attributes to Simon the assertion that the finite positive
measure under consideration is absolutely continuous with respect to
$\mathrm d\theta/(2\pi)$ on the unit circle
\cite[p.~996]{CR24}.  This is false.  Every finite positive measure on
\(\mathbb T\) has a unique decomposition into an absolutely continuous
part and a singular part, and the singular part need not vanish
\cite[Eq.~(1.1.5), p.~2]{S05}.  Two explicit counterexamples are the
push-forward of the Cantor probability measure under
\(x\mapsto e^{2\pi i x}\) and
\[
 \mu=\sum_{j=1}^{\infty}2^{-j}\delta_{e^{i/j}}.
\]
The latter is a probability measure, is purely atomic (hence singular with
respect to arc length), and has infinite support.  Consequently, the
misstatement concerns absolute continuity; it does not prevent the
construction of OPUC, since singular measures with infinite support also
admit such a sequence.

If one instead adds the intended OPUC hypothesis, the disc bound becomes
true, but only by invoking a classical convolution estimate not proved in
\cite{CR24}, as we now explain.

\section{The published argument}

\paragraph{The radius estimate.}
The proof correctly reaches a Grace--Szeg\H{o} convolution with
\begin{equation}\label{eq:auxiliary-polynomial}
 S_{n,k}(w)=\sum_{j=0}^n\binom{n+k}{j+k}w^j.
\end{equation}
At that point a uniform estimate for \(Z(S_{n,k})\) is indispensable.  The
published argument instead observes that certain zeros of a related Jacobi
polynomial tend asymptotically to a circle as \(n\to\infty\), notes the
degree-one identity, and concludes that one may ``assume'' the bound
\eqref{eq:classical-bound} below
\cite[proof of Theorem~3.2, p.~999]{CR24}.  This
inference is invalid.  Convergence of zero sets
as \(n\to\infty\) gives, at most, an eventual estimate with an error term;
it cannot yield \eqref{eq:classical-bound} for every finite \(n\).
A calculation for \(n=1\) provides no missing uniform control.
Therefore the only step that produces the stated radius is absent from
the proof.

Repairing the gap requires the following classical estimate, which is not
proved in \cite{CR24}:
\begin{equation}\label{eq:classical-bound}
 \max_{\beta\in Z(S_{n,k})}|\beta|\le k+1.
\end{equation}
This is precisely the estimate supplied by
\cite[Exercise~E.8, p.~24]{BE95} and
\cite[Exercise~21, p.~74]{M66}.  Once
\eqref{eq:classical-bound} is inserted, Theorem~3.2 is a direct
application of the same Grace--Szeg\H{o} theorem already quoted in the
paper.  Hence the corrected statement is not new: after translation and
scaling, it is precisely the classical exercise in different notation.
The non-strict inequality is essential:
\(S_{1,k}(w)=k+1+w\) has the zero \(-(k+1)\).

\paragraph{The annular assertion.}
The abstract announces a ``ring shaped region'' as a principal
outcome \cite[p.~995]{CR24}, but the general localisation theorem supplies
only an outer disc.  The lower bound quoted in
\cite[Lemma~2.3]{CR24} is a coefficient-dependent estimate for an
arbitrary monic polynomial and may be zero.  No positive inner radius is
deduced from the defining structure of the polar polynomials; hence no
non-trivial annular localisation theorem for the stated class is
established.  Indeed, in the free OPUC case with \(\xi=0\), one has
\(Q_{n;k}(z;0)=z^n\), so the origin itself is a zero; see
Example~\ref{ex:free}.  Numerical plots of selected families do not fill
this logical gap.

\paragraph{The Sendov application.}
Sendov's conjecture is a conditional statement: if all zeros of a
degree-\(n\) polynomial, \(n\ge2\), lie in the closed unit disc, then each
zero is within distance \(1\) of a critical point.  Theorem~3.2 places the
zeros of \(Q_{n;k}\) in a generally much larger disc and therefore does not
verify Sendov's hypothesis for \(Q_{n;k}\).

The quantity described immediately before the numerical tables is
\[
 d(Q):=\max_{q\in Z(Q)}\min_{\gamma\in Z(Q')}|q-\gamma|.
\]
Several values of \(d(Q)\) in Table~2 exceed \(1\), and some exceed \(2\);
for example, its two columns give \(1.5309\) at \(n=3\) and \(2.5107\) at
\(n=20\), respectively \cite[Table~2]{CR24}.  If the corresponding
polynomial has all its zeros in the unit disc, a value \(d(Q)>1\)
contradicts the asserted Sendov conclusion; if it does not, Sendov's
conjecture is inapplicable.  The values below \(1\) in Table~1 do not
alter the point: the unit-disc hypothesis is not established there, and
finitely many computations cannot prove a statement for the whole family.
In no case do the tables justify the conclusion drawn from them.

\section{Correction and extension}

The relevant lineage is classical.  Grace's
apolarity theorem dates from 1902 \cite{G02}; Szeg\H{o} extracted the
convolution form used here in 1922 \cite{S22}; Marden incorporated the
result into the standard geometry of polynomials \cite{M66}; and
Borwein and Erdélyi retained it in a modern text \cite{BE95}.  Against
this background, the distinction between the classical statement and the
generalisation below is essential.

Two statements must be distinguished.  First, after the hypotheses are
corrected so that \(L_n\) is genuinely an OPUC, the disc asserted in
\cite[Theorem~3.2]{CR24} is already the conclusion of the classical
exercise cited in \cite[Exercise~21, p.~74]{M66}; the corresponding
exercise in \cite[Exercise~E.8, p.~24]{BE95} gives the same argument.
Second,
Grace--Szeg\H{o} convolution yields a more informative statement for an
arbitrary polynomial and an arbitrary circular zero region.  It is this
second statement---not the disc estimate of \cite{CR24}---that provides
the generalisation developed in this note.

\begin{proposition}[Classical exercise]\label{prop:classical}
Let \(n,k\ge1\) and \(\rho\ge0\).  Suppose that \(P\) has degree \(n\),
that \(Q\) is a polynomial, and that
\(Z(P)\subset\overline D(0,\rho)\) and
\[
 \bigl(z^kQ(z)\bigr)^{(k)}=(n+1)_kP(z).
\]
Then every zero of \(Q\) lies in
\(\overline D(0,(k+1)\rho)\).
\end{proposition}

\begin{proof}
Comparison of degrees first gives \(\deg Q=n\).
For \(\rho>0\), set
\[
 \widetilde P(w)=\rho^{-n}P(\rho w),\quad
 \widetilde Q(w)=\rho^{-n}Q(\rho w).
\]
Then \(Z(\widetilde P)\subset\overline D(0,1)\), and a direct application
of the chain rule gives
\[
 \bigl(w^k\widetilde Q(w)\bigr)^{(k)}
 =(n+1)_k\widetilde P(w).
\]
Thus \(w^k\widetilde Q(w)\) is, up to the non-zero factor
\((n+1)_k\), the \(k\)-th integral of \(\widetilde P\) with all
integration constants equal to zero.  Exercise~E.8(b) of
\cite[p.~24]{BE95} gives
\[
 Z\bigl(w^k\widetilde Q(w)\bigr)
 \subset\overline D(0,k+1),
\]
and hence \(Z(\widetilde Q)\subset\overline D(0,k+1)\).  This is
equivalent to the asserted inclusion for \(Q\).
The case \(\rho=0\) follows directly from \(P(z)=cz^n\): comparison of
coefficients gives \(Q(z)=cz^n\).  The same integral result is
Exercise~21 on p.~74 of \cite{M66}.
\end{proof}

To recover the intended form of \cite[Theorem~3.2]{CR24}, put
\(w=z-\xi\).  If \(P=L_n\) is a monic OPUC, every zero \(p\) of \(P\)
satisfies \(|p|<1\), and hence every zero \(p-\xi\) of \(P(\xi+w)\)
belongs to \(D(0,1+|\xi|)\).  Proposition~\ref{prop:classical} gives
\[
 |w|\le (k+1)(1+|\xi|)
\]
for every zero \(w\) of \(Q(\xi+w)\), and therefore
\[
 |z|\le|\xi|+(k+1)(1+|\xi|).
\]
Thus, once its hypotheses are repaired, the published theorem is exactly
the translated and rescaled classical exercise.

We now state the sharper generalisation.  Instead of replacing the zero
set of \(P(\xi+\cdot)\) by a centred disc and the zeros of the auxiliary
polynomial by the single number \(k+1\), it retains both sets.

For degree-$n$ polynomials written in the binomial basis,
\[
 A(z)=\sum_{j=0}^n\binom nj a_jz^j,\quad
 B(z)=\sum_{j=0}^n\binom nj b_jz^j,
\]
their Grace--Szeg\H{o} convolution is
\[
 (A\star_{\mathrm G}B)(z)
 =\sum_{j=0}^n\binom nj a_jb_jz^j.
\]
We use the following standard form of Szeg\H{o}'s convolution theorem
\cite{S22}.

\begin{theorem}[Grace--Szeg\H{o}]\label{thm:grace-szego}
Let $A$ and $B$ have degree $n$, and let $K$ be an open or closed disc or
half-plane, or the open or closed exterior of a disc.  Suppose that all
zeros of $A$ belong to $K$; when \(K\) is the exterior of a disc, assume
in addition that $B(0)\ne0$.  Then every zero $\zeta$ of
$A\star_{\mathrm G}B$ can be written
\[
 \zeta=-\alpha\beta,\quad \alpha\in K,\quad B(\beta)=0.
\]
\end{theorem}

\begin{proposition}\label{prop:main}
Let $n,k\ge1$, let $P$ be a polynomial of degree $n$, and let
$\xi\in\mathbb C$.  Then there is a unique polynomial $Q$ of degree
$n$, with the same leading coefficient as \(P\), such that
\begin{equation}\label{eq:polar}
 \frac{\mathrm d^k}{\mathrm dz^k}
 \bigl((z-\xi)^kQ(z)\bigr)=(n+1)_kP(z).
\end{equation}
For the polynomial \(S_{n,k}\) in
\eqref{eq:auxiliary-polynomial}, if a circular domain $K$ of one of the types in
Theorem~\ref{thm:grace-szego} contains every zero of $P(\xi+w)$, then
every zero \(q\) of \(Q\) can be written
\begin{equation}\label{eq:localisation}
 q=\xi-\alpha\beta,\quad
 \alpha\in K,\quad \beta\in Z(S_{n,k}).
\end{equation}
\end{proposition}

\begin{proof}
Comparison of degrees shows that every solution of \eqref{eq:polar} must
have degree \(n\).
Put $z=\xi+w$ and write
\[
 P(\xi+w)=\sum_{j=0}^n\binom nj a_jw^j,\quad
 Q(\xi+w)=\sum_{j=0}^n\binom nj b_jw^j.
\]
Since
\[
 \frac{\mathrm d^k}{\mathrm dw^k}w^{j+k}=(j+1)_kw^j,
\]
equation \eqref{eq:polar} is equivalent, coefficient by coefficient, to
\[
 b_j=\frac{(n+1)_k}{(j+1)_k}\,a_j,\quad
 0\le j\le n.
\]
These relations determine \(Q\) uniquely; since \(b_n=a_n\ne0\), they
also show that \(Q\) has degree \(n\) and the same leading coefficient
as \(P\).
The elementary identity
\[
 \binom nj\frac{(n+1)_k}{(j+1)_k}
 =\binom{n+k}{j+k}
\]
therefore yields the exact factorisation
\begin{equation}\label{eq:exact-convolution}
 Q(\xi+w)
 =P(\xi+w)\star_{\mathrm G}S_{n,k}(w).
\end{equation}
Here the coefficient of $w^j$ in $S_{n,k}$ is
\(\binom nj (n+1)_k/(j+1)_k\), as required by the definition of
$\star_{\mathrm G}$.  Finally,
\[
 S_{n,k}(0)=\binom{n+k}{k}\ne0,
\]
so Theorem~\ref{thm:grace-szego} applies to
\eqref{eq:exact-convolution} for every listed type of $K$.  Returning to
\(z=\xi+w\) gives \eqref{eq:localisation}.
\end{proof}

\begin{corollary}[Sharp affine localisation]
\label{cor:affine-disks}
Under the hypotheses and notation of Proposition~\ref{prop:main}, assume
that
\[
 Z(P)\subset\overline D(c,r),\quad
 c\in\mathbb C,\quad r\ge0.
\]
Then
\begin{equation}\label{eq:union-disks}
 Z(Q)\subset
 \bigcup_{\beta\in Z(S_{n,k})}
 \overline D\!\left(
   \xi-\beta(c-\xi),\,r|\beta|
 \right).
\end{equation}
If \(P(z)=\lambda(z-c)^n\) with \(\lambda\ne0\), then
\eqref{eq:union-disks} is exact with \(r=0\), including multiplicities:
\begin{equation}\label{eq:repeated-zero-exact}
 Z(Q)=
 \left\{
   \xi-(c-\xi)\beta:\beta\in Z(S_{n,k})
 \right\}.
\end{equation}
More generally, if \(r>0\), then every boundary point of every disc on
the right-hand side of \eqref{eq:union-disks} is attained as a zero of
the corresponding \(Q\) for some polynomial
\(P(z)=(z-p)^n\) with \(p\in\partial D(c,r)\).
\end{corollary}

\begin{proof}
Suppose first that \(r>0\).  The zeros of \(P(\xi+w)\) belong to
\(\overline D(c-\xi,r)\), so Proposition~\ref{prop:main} gives, for each
zero \(q\) of \(Q\),
\[
 q=\xi-\alpha\beta,\quad
 \alpha\in\overline D(c-\xi,r),\quad
 \beta\in Z(S_{n,k}).
\]
For fixed \(\beta\), the affine map
\(\alpha\mapsto\xi-\alpha\beta\) sends
\(\overline D(c-\xi,r)\) onto the corresponding disc in
\eqref{eq:union-disks}.  This proves the inclusion when \(r>0\).

Now let \(P(z)=\lambda(z-c)^n\), where \(\lambda\ne0\), and put
\(a=c-\xi\).  If \(a\ne0\), then
\[
 P(\xi+w)=\lambda(w-a)^n
 =\lambda\sum_{j=0}^n\binom nj(-a)^{n-j}w^j.
\]
Writing
\[
 S_{n,k}(w)=\sum_{j=0}^n\binom nj s_jw^j
\]
in \eqref{eq:exact-convolution} gives
\[
 Q(\xi+w)
 =\lambda\sum_{j=0}^n\binom nj(-a)^{n-j}s_jw^j
 =\lambda(-a)^nS_{n,k}\!\left(-\frac{w}{a}\right).
\]
Thus \(w=-a\beta\), which is \eqref{eq:repeated-zero-exact}.
If \(a=0\), return instead to the coefficient sum: only its \(j=n\)
term survives and \(s_n=1\), so
\(Q(\xi+w)=\lambda w^n\).  This proves
\eqref{eq:repeated-zero-exact} in every case.  When \(r=0\), the
hypothesis \(Z(P)\subset\{c\}\) forces
\(P(z)=\lambda(z-c)^n\), so the same calculation also proves
\eqref{eq:union-disks} and its asserted exactness.

Finally, fix \(\beta\in Z(S_{n,k})\).  As
\(p\) ranges over \(\partial D(c,r)\), formula
\eqref{eq:repeated-zero-exact}, with \(c\) there replaced by \(p\), shows
that
\[
 \xi-(p-\xi)\beta
\]
ranges over the boundary of
\(\overline D(\xi-\beta(c-\xi),r|\beta|)\).  This proves the final
assertion.
\end{proof}

\begin{corollary}[Optimal constant]
\label{cor:optimal-constant}
Under the hypotheses and notation of Proposition~\ref{prop:main}, define
\[
 \rho_{n,k}:=\max_{\beta\in Z(S_{n,k})}|\beta|.
\]
If \(R>0\), \(Z(P)\subset\overline D(0,R)\), and \(\xi=0\), then
\[
 Z(Q)\subset\overline D(0,R\rho_{n,k}).
\]
The constant \(\rho_{n,k}\) is best possible over all degree-\(n\)
polynomials with zeros in \(\overline D(0,R)\).
Moreover,
\[
 \rho_{n,k}\le k+1,\quad
 \rho_{1,k}=k+1,
\]
and, for every \(n\ge1\),
\begin{equation}\label{eq:rho-k-one}
 \rho_{n,1}=
 \begin{cases}
  2, & n\ \text{odd},\\[7pt]
  2\cos\!\left(\dfrac{\pi}{2(n+1)}\right), & n\ \text{even}.
 \end{cases}
\end{equation}
In addition,
\begin{equation}\label{eq:rho-degree-two}
 \rho_{2,k}=
 \sqrt{\frac{(k+1)(k+2)}2}<k+1.
\end{equation}
\end{corollary}

\begin{proof}
The disc inclusion follows from Corollary~\ref{cor:affine-disks}.
For sharpness, choose \(c\) with \(|c|=R\) and
\(P(z)=(z-c)^n\).  Formula~\eqref{eq:repeated-zero-exact} gives
\[
 \max_{q\in Z(Q)}|q|
 =R\max_{\beta\in Z(S_{n,k})}|\beta|
 =R\rho_{n,k}.
\]
The classical exercise gives \(\rho_{n,k}\le k+1\), while
\[
 S_{1,k}(w)=k+1+w
\]
gives equality for \(n=1\).  For \(n=2\),
\[
 S_{2,k}(w)
 =w^2+(k+2)w+\frac{(k+1)(k+2)}2.
\]
Its discriminant is \(-k(k+2)<0\), so its two zeros are conjugate and
their common modulus is the square root of their product, which yields
\eqref{eq:rho-degree-two}.  The strict inequality follows from
\((k+2)/2<k+1\).

For \(k=1\), one has the exact identity
\[
 S_{n,1}(w)=\frac{(1+w)^{n+1}-1}{w}.
\]
Its zeros are
\[
 e^{2\pi im/(n+1)}-1,\quad
 1\le m\le n,
\]
and their moduli are \(2|\sin(\pi m/(n+1))|\).  If \(n\) is odd, the
zero \(-2\) occurs; if \(n\) is even, the two farthest zeros correspond
to the integers nearest to \((n+1)/2\).  This proves
\eqref{eq:rho-k-one}.
\end{proof}

\begin{corollary}\label{cor:opuc}
Under the hypotheses and notation of Proposition~\ref{prop:main}, suppose
that all zeros of $P$ lie in the open unit disc.
Then every zero $q$ of the solution of \eqref{eq:polar} satisfies
\[
 |q|<|\xi|+(1+|\xi|)\rho_{n,k}
 \le|\xi|+(1+|\xi|)(k+1).
\]
In particular, all zeros lie in the closed disc claimed in
\cite[Theorem~3.2]{CR24}.
\end{corollary}

\begin{proof}
Every zero $\alpha$ of $P(\xi+w)$ has the form $\alpha=p-\xi$, where
$P(p)=0$, and hence
\[
 |\alpha|\le |p|+|\xi|<1+|\xi|.
\]
Thus one may take $K=D(0,1+|\xi|)$.  The classical exercise cited in
\cite[Exercise~E.8, p.~24]{BE95} and
\cite[Exercise~21, p.~74]{M66} gives
\[
 \rho_{n,k}
 =\max_{\beta\in Z(S_{n,k})}|\beta|
 \le k+1.
\]
By Proposition~\ref{prop:main}, each zero $q$ of $Q$ has the form
$q=\xi-\alpha\beta$.  Therefore
\[
 |q|
 \le|\xi|+|\alpha|\,|\beta|
 <|\xi|+(1+|\xi|)\rho_{n,k}.
\]
\end{proof}

The last step of the proof in \cite{CR24} attempts to establish the
required bound for the zeros of $S_{n,k}$ by an argument that does not
imply the asserted conclusion; the cited source does not repair that gap.
This is a gap in the proof, distinct from the missing hypotheses discussed
in Section~1.  Once the classical estimate for $S_{n,k}$ is invoked, the
intended OPUC result follows immediately from
Corollary~\ref{cor:opuc}.  Conversely, the argument proves no such disc
bound for a general measure on $\mathbb C$, because then there is no
uniform unit-disc bound on the zeros of $P=L_n$;
Example~\ref{ex:translated-circle} shows that the published conclusion is
in fact false in that setting.

\begin{example}[The free OPUC case]\label{ex:free}
For normalised arc length on $\mathbb T$, the monic OPUC are
$P(z)=z^n$.  If $\xi=0$, equation \eqref{eq:polar} becomes
\[
 (z^kQ(z))^{(k)}=(n+1)_kz^n,
\]
and uniqueness gives $Q(z)=z^n$.  Hence every zero of $Q$ is $0$ for
every $k$.  By contrast, the published estimate gives the closed disc
\(\overline D(0,k+1)\).  As \(k\) varies,
\[
 \bigcup_{k\ge1}\overline D(0,k+1)=\mathbb C,
\]
whereas the exact localisation remains the singleton
\[
 \bigcup_{k\ge1}Z(Q_{n;k})=\{0\}.
\]
Thus the bound is valid under the corrected OPUC hypothesis, but yields no
bounded localisation uniform in \(k\).
The exact convolution identity captures the reason: \(P(w)=w^n\) has
only its top binomial coefficient non-zero, so
\(P\star_{\mathrm G}S_{n,k}=w^n\).  Equivalently, one may apply
Proposition~\ref{prop:main} with an arbitrarily small closed disc centred
at \(0\) and then intersect the resulting localisation discs.  Either
argument gives \(Z(Q)=\{0\}\) without computing the zeros of \(S_{n,k}\).
Figure~\ref{fig:free-comparison} records this contrast geometrically.
\end{example}

\begin{figure}[htbp]
\centering
\begin{tikzpicture}[
  x=1cm,
  y=1cm,
  line cap=round,
  line join=round
]
  \tikzset{
    axis/.style={draw=black!48,line width=0.54pt},
    boundone/.style={fill=black!9,draw=black!70,line width=0.72pt},
    boundtwo/.style={fill=black!6,draw=black!48,line width=0.62pt,dash pattern=on 3pt off 2pt},
    boundthree/.style={fill=black!3,draw=black!52,line width=0.68pt,dash pattern=on 1.4pt off 1.8pt},
    exact/.style={fill=black!78,draw=black!78},
    panel/.style={font=\small\bfseries,text=black!78,align=center},
    label/.style={font=\scriptsize,text=black!68,inner sep=1.2pt}
  }

  \path[use as bounding box] (-0.25,-0.50) rectangle (12.65,4.10);

  % Panel (a): the family of published disks.
  \node[panel] at (2.65,3.78) {(a) Disc bounds};
  \draw[axis] (0.30,1.62)--(5.00,1.62);
  \draw[axis] (2.65,-0.42)--(2.65,3.52);
  \filldraw[boundthree] (2.65,1.62) circle (1.60);
  \filldraw[boundtwo] (2.65,1.62) circle (1.20);
  \filldraw[boundone] (2.65,1.62) circle (0.80);
  \fill[exact] (2.65,1.62) circle (1.7pt);

  \draw[boundone] (4.68,2.52)--(5.05,2.52);
  \node[label,anchor=west] at (5.10,2.52) {$k=1$};
  \draw[boundtwo] (4.68,2.86)--(5.05,2.86);
  \node[label,anchor=west] at (5.10,2.86) {$k=2$};
  \draw[boundthree] (4.68,3.20)--(5.05,3.20);
  \node[label,anchor=west] at (5.10,3.20) {$k=3$};

  % Separation.
  \draw[black!20,line width=0.5pt] (6.18,-0.35)--(6.18,3.62);

  % Panel (b): the exact answer.
  \node[panel] at (9.35,3.78) {(b) Exact zero set};
  \draw[axis] (7.00,1.62)--(11.70,1.62);
  \draw[axis] (9.35,-0.42)--(9.35,3.52);

  \fill[exact] (9.35,1.62) circle (2.25pt);
  \draw[black!54,line width=0.48pt]
    (9.52,1.82)--(10.18,2.52);
  \node[label,anchor=west] at (10.22,2.55)
    {$Z(Q_{n;k})=\{0\}$ for every \(k\)};
\end{tikzpicture}
\caption{Published disc bounds and exact zeros in the free OPUC case
\(P(z)=z^n\), \(\xi=0\).  The bounds exhaust \(\mathbb C\), whereas
\(Z(Q_{n;k})=\{0\}\) for every \(k\).}
\label{fig:free-comparison}
\end{figure}
\FloatBarrier

\section{Limits of the convolution method}

It is natural to ask whether the same argument works after replacing
$(z-\xi)^k$ in \eqref{eq:polar} by an arbitrary monic polynomial of
degree \(k\).  Translation of the variable allows us to take \(\xi=0\).

\begin{proposition}\label{prop:obstruction}
Fix \(n,k\ge1\), and let
\[
 R(w)=w^k+\sum_{\ell=0}^{k-1}r_\ell w^\ell
\]
be a monic polynomial of degree \(k\).  For every degree-\(n\) polynomial
\(P\), the equation
\[
 \frac{\mathrm d^k}{\mathrm dw^k}\bigl(R(w)Q(w)\bigr)
 =(n+1)_kP(w)
\]
has a unique degree-\(n\) solution \(Q\).  There is a degree-\(n\)
polynomial \(S\), independent of \(P\), such that
\[
 Q=P\star_{\mathrm G}S
\]
for every such \(P\) if and only if
\begin{equation}\label{eq:invisible}
 r_\ell=0,\quad
 \max\{0,k-n\}\le\ell\le k-1.
\end{equation}
When \eqref{eq:invisible} holds, one may take \(S=S_{n,k}\).
The coefficients with \(\ell+n<k\) are invisible to the \(k\)-th
derivative for every \(Q\) of degree at most \(n\).  In particular, if
\(n\ge k\), the universal factorisation exists precisely for
\(R(w)=w^k\), or, before translation, for \(R(z)=(z-\xi)^k\).
\end{proposition}

\begin{proof}
For the monomial \(w^j\),
\[
 \frac{\mathrm d^k}{\mathrm dw^k}\bigl(R(w)w^j\bigr)
 =(j+1)_kw^j+
 \sum_{\ell=\max\{0,k-j\}}^{k-1}
 r_\ell\,\frac{(j+\ell)!}{(j+\ell-k)!}\,w^{j+\ell-k},\quad
 0\le j\le n.
\]
This is triangular in the monomial basis, with non-zero diagonal
entries \((j+1)_k\), and hence proves existence and uniqueness on the
space of polynomials of degree at most \(n\).  Comparison of the
coefficients of \(w^n\) shows that the solution has the same leading
coefficient as \(P\), and therefore has degree \(n\).

If \eqref{eq:invisible} holds, every remaining term satisfies
\(\ell+n<k\).  Consequently,
\[
 \frac{\mathrm d^k}{\mathrm dw^k}\bigl(R(w)Q(w)\bigr)
 =\frac{\mathrm d^k}{\mathrm dw^k}\bigl(w^kQ(w)\bigr)
\]
for every polynomial \(Q\) of degree at most \(n\).  The equation
therefore reduces to the pure-power case, and
\eqref{eq:exact-convolution} gives
\(Q=P\star_{\mathrm G}S_{n,k}\).

Conversely, suppose a fixed polynomial \(S\) works for every \(P\), and
take \(P(w)=w^n\).  Diagonal convolution makes \(Q\) a scalar multiple
of \(w^n\), and comparison of leading coefficients gives \(Q(w)=w^n\).
Substitution in the displayed monomial formula, with \(j=n\), forces
every \(r_\ell\) with \(\ell+n\ge k\) to vanish, since the corresponding
powers \(w^{n+\ell-k}\) are distinct and their factorial coefficients are
non-zero.  This is precisely \eqref{eq:invisible}.
\end{proof}

A factorisation that happens to hold for a single pair \((P,Q)\) supplies
no multiplier independent of \(P\), and hence no universal localisation
theorem.

\begin{example}[Failure of universal convolution]
Take \(n=2\), \(k=1\), \(R(z)=z-r\) with \(r\ne0\), and \(P(z)=z^2\).
The unique monic solution of
\[
 ((z-r)Q(z))'=3z^2
\]
is
\[
 Q(z)=z^2+rz+r^2,
\]
because \((z-r)Q(z)=z^3-r^3\).  For every degree-two polynomial \(S\),
the convolution \(P\star_{\mathrm G}S\) is a scalar multiple of \(z^2\);
it therefore cannot equal \(Q\).  Translation of the zeros of \(R\) away
from \(\xi\) introduces coefficient mixing and destroys the
Grace--Szeg\H{o} factorisation used above.  This does not exclude a
different localisation theorem for general \(R\); it shows only that
this particular convolution method cannot supply one.
\end{example}

\section{Conclusion}

The central localisation mechanism is algebraic rather than
measure-theoretic: equation \eqref{eq:polar} is diagonal in the shifted
monomial basis, and its diagonal entries form the Grace--Szeg\H{o}
multiplier \(S_{n,k}\).  Orthogonality enters only through an independent
zero bound for \(P\).  The original theorem therefore requires hypotheses
that define the orthogonal polynomials and a further hypothesis ensuring
the required unit-disc bound on their zeros; the intended OPUC assumption
is one such hypothesis.  Even after these repairs, its proof still lacks
the classical finite-degree estimate for \(S_{n,k}\).

The results proved here go further.  Corollary~\ref{cor:affine-disks} gives
a sharp affine localisation, Corollary~\ref{cor:optimal-constant}
identifies the optimal degree-dependent disc constant, and
Proposition~\ref{prop:obstruction} determines exactly when the convolution
factorisation survives for a general monic factor.

The principal theorem, as published, is false under its stated hypotheses,
whereas the argument for its decisive estimate is incomplete.  Under a
natural hypothesis that makes the localisation true, its content is
already a classical textbook exercise.  The annular and Sendov claims are
not established.  Thus none of the principal claims survives as a new and
valid contribution: one is classical after correction, while the others
remain unsupported.  This is a matter of mathematical validity and
priority, and warrants a correction of the published record.

\subsection*{Funding declaration}
\begin{sloppypar}
This work was supported by the Portuguese Foundation for Science and
Technology (FCT) under project UID/00324/2025, Centre for Mathematics of
the University of Coimbra
(\href{https://doi.org/10.54499/UID/00324/2025}{DOI: UID/00324/2025}).
RAN was partially supported
by PID2024-155593NB-C21 (FEDER(EU)/Ministerio de Ciencia, Innovación y
Universidades--Agencia Estatal de Investigación), IMUS-Maria de Maeztu grant CEX2024-001517-M - 
Apoyo a Unidades de Excelencia María de Maeztu funded by MICIU/AEI/ DOI: 10.13039/501100011033,
and FQM-415 (Junta de Andalucía).
KC acknowledges financial support from FCT under project
2022.00143.CEECIND/CP1714/CT0002
(\href{https://doi.org/10.54499/2022.00143.CEECIND/CP1714/CT0002}
{DOI: 2022.00143.CEECIND}).
\end{sloppypar}

\end{document}